\documentclass{amsart}
\usepackage{amsmath} 
\usepackage{amssymb}
\usepackage{mathtools}
\usepackage{centernot}
\usepackage{stmaryrd}
\usepackage{esvect}
\usepackage{amsthm}
\usepackage{upgreek}
\usepackage{comment}
\usepackage{amsfonts}
\usepackage{physics}
\usepackage{tikz-cd}
\usepackage[utf8]{inputenc}
\usepackage[english]{babel}
\usepackage{graphicx}
\usepackage[margin=1in]{geometry}
\usepackage{mathtools}
\usepackage{hyperref}
\usepackage{microtype}

\graphicspath{{images/}}
\newcommand{\R}{\mathbb{R}}

\renewcommand{\part}[1]{\noindent\textit{Part #1)}}

\newcommand{\ba}{\begin{align*}}

\newcommand{\ea}{\end{align*}}

\newcommand{\N}{\mathbb{N}}

\newcommand{\Int}[1]{%
  {\kern0pt#1}^{\mathrm{o}}%
}

\newcommand{\Z}{\mathbb{Z}}

\newcommand{\bp}{\begin{pmatrix}}
\newcommand{\ep}{\end{pmatrix}}

\newcommand{\interior}[1]{%
  {\kern0pt#1}^{\mathrm{o}}%
}

\usepackage{xargs}                      
\usepackage{xcolor}
\usepackage[colorinlistoftodos,prependcaption,textsize=tiny]{todonotes}
\newtheorem{theorem}{Theorem}[section]
\newtheorem{lemma}[theorem]{Lemma}

\theoremstyle{definition}
\newtheorem{definition}[theorem]{Definition}
\newtheorem{algorithm}[theorem]{Algorithm}
\newtheorem{conjecture}[theorem]{Conjecture}
\newtheorem{corollary}[theorem]{Corollary}
\newtheorem{question}[theorem]{Question}
\newtheorem{innercustomthm}{Theorem}
\newenvironment{customthm}[1]
  {\renewcommand\theinnercustomthm{#1}\innercustomthm}
  {\endinnercustomthm}

\title{Finiteness of Singular Orbits of Pseudo Anosov Flows}

\author[L. Li]{Lily Qiao Li}
\address{Department of Mathematics, Princeton University, Fine Hall, 304 Washington Rd, Princeton, NJ 08544}
\email{lilyli@princeton.edu}

\subjclass[2020]{}
\keywords{}

\begin{document}  

\begin{abstract} We show that for any closed atoroidal 3-manifold $M$, there are only finitely many isotopy classes of links that can arise as singular orbits of pseudo-Anosov flows on $M$.
\end{abstract}

\maketitle

\section{Introduction}

In this paper, we characterize links that arise as singular orbits of pseudo-Anosov flows on an atoroidal 3-manifold $M$. In particular, we show that the set of isotopy classes of such links is finite for any fixed $M$. We also show combinatorial finiteness for all possible gut regions of stable laminations of pseudo-Anosov flows on a fixed $M$. The proofs have a combinatorial flavor and draws from techniques in essential lamination theory and normal surface theory.

Pseudo-Anosov flows are flows on 3-manifolds that behave like Anosov flows away from finitely many singular orbits. They were originally introduced by Thurston as a natural generalization of Anosov flows. Informally, they are locally modeled on suspension flows of pseudo-Anosov homeomorphisms on a surface of genus $g$.  

Pseudo-Anosov flows have been shown to be abundant. In particular, for any $n\in\N$, there are 3-manifolds with $n$ distinct pseudo-Anosov flows up to orbit equivalence \cite{BBYabundance}\cite{BMabundance}\cite{CPabundance}\cite{BYabundance}\cite{BSZabundance}. This paper is motivated by the following finiteness conjecture for Anosov and pseudo-Anosov flows. 

\begin{conjecture}[The Finiteness Conjecture]
    There are only finitely many transitive Anosov and pseudo-Anosov flows up to orbit equivalence on any closed 3-manifold.
\end{conjecture}

There has been much activity and progress over the years on this conjecture from various viewpoints. More recently, work of Barthelm\'e--Bowden--Mann proved the finiteness conjecture for Reeb Anosov flows using the classification of transitive pseudo-Anosov flows up to orbit equivalence via free homotopy classes of closed orbits due to Barthelm\'e--Mann and Barthelm\'e--Frankel--Mann \cite{BM24}\cite{BFMclassification}. Zung \cite{Zungfiniteness} and Baldwin--Sivek--Zung \cite{BSZabundance} extended this to show finiteness for pseudo-Anosov flows admitting positive Birkhoff sections in rational homology spheres. Chaidez--Pan \cite{ChaidezPan} proved the finiteness conjecture for pseudo-Anosov Reeb flows, which they define. 

In this paper, we provide evidence to support the finiteness conjecture, characterizing the links of singular orbits that arise for a fixed atoroidal 3-manifold $M$. 

\begin{customthm}{A}\label{main}
    Let $M$ be a closed atoroidal manifold. Then the set of isotopy classes of links representing singular orbits of pseudo-Anosov flows on M is finite. Further, the set of isotopy classes of degeneracy curves of these pseudo-Anosov flows is also finite.
\end{customthm}

Barthelm\'e–Tsang–Zung \cite{noperfectfits} have recently used Theorem A as the first step in proving the finiteness conjecture for pseudo-Anosov flows without perfect fits. Other works in progress on this conjecture include the following: Landry--Taylor \cite{LandryTaylor} show that there are only finitely many almost transverse pseudo-Anosov flows for every finite-depth foliation in a hyperbolic 3-manifold, and work of Barthelm\'e--Paulet \cite{BarthelmePaulet} proves the finiteness conjecture for transitive pseudo-Anosov flows in graph-manifolds up to finite covers.

\begin{figure}
    \centering
    \includegraphics[scale = 0.11]{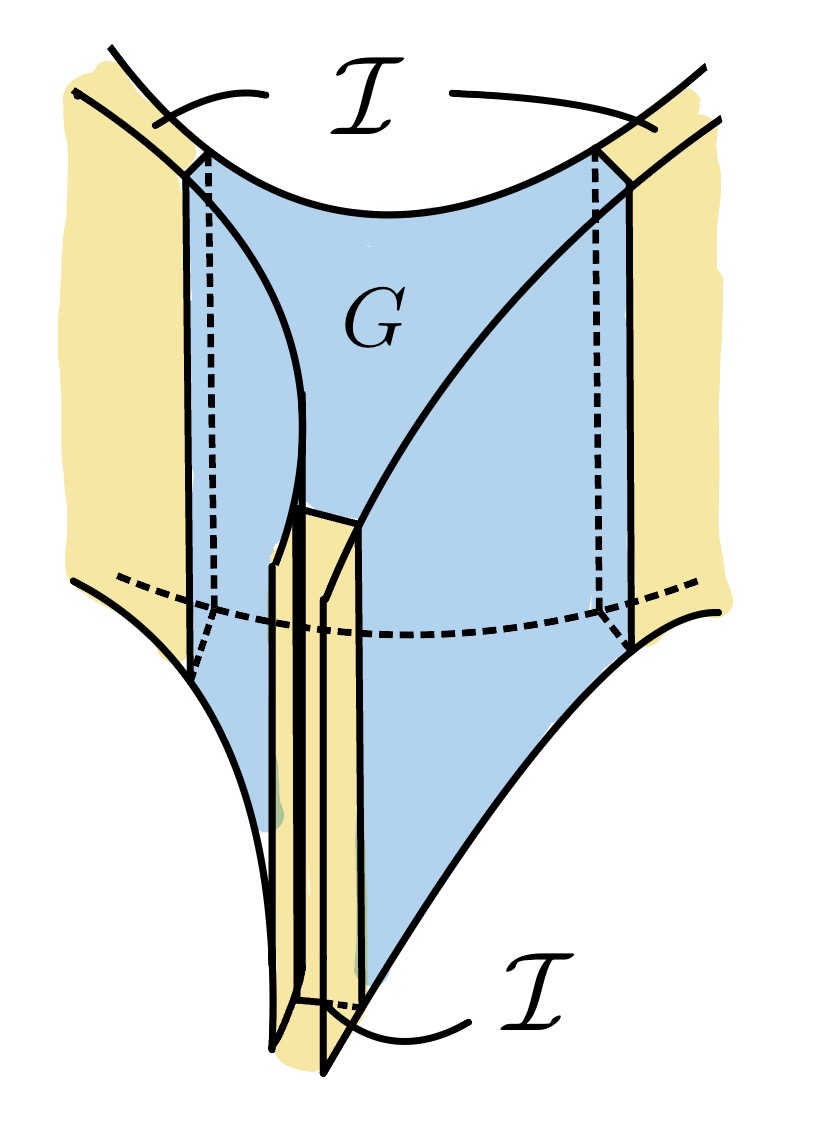}
    \hspace{3em}
    \includegraphics[scale = 0.09]{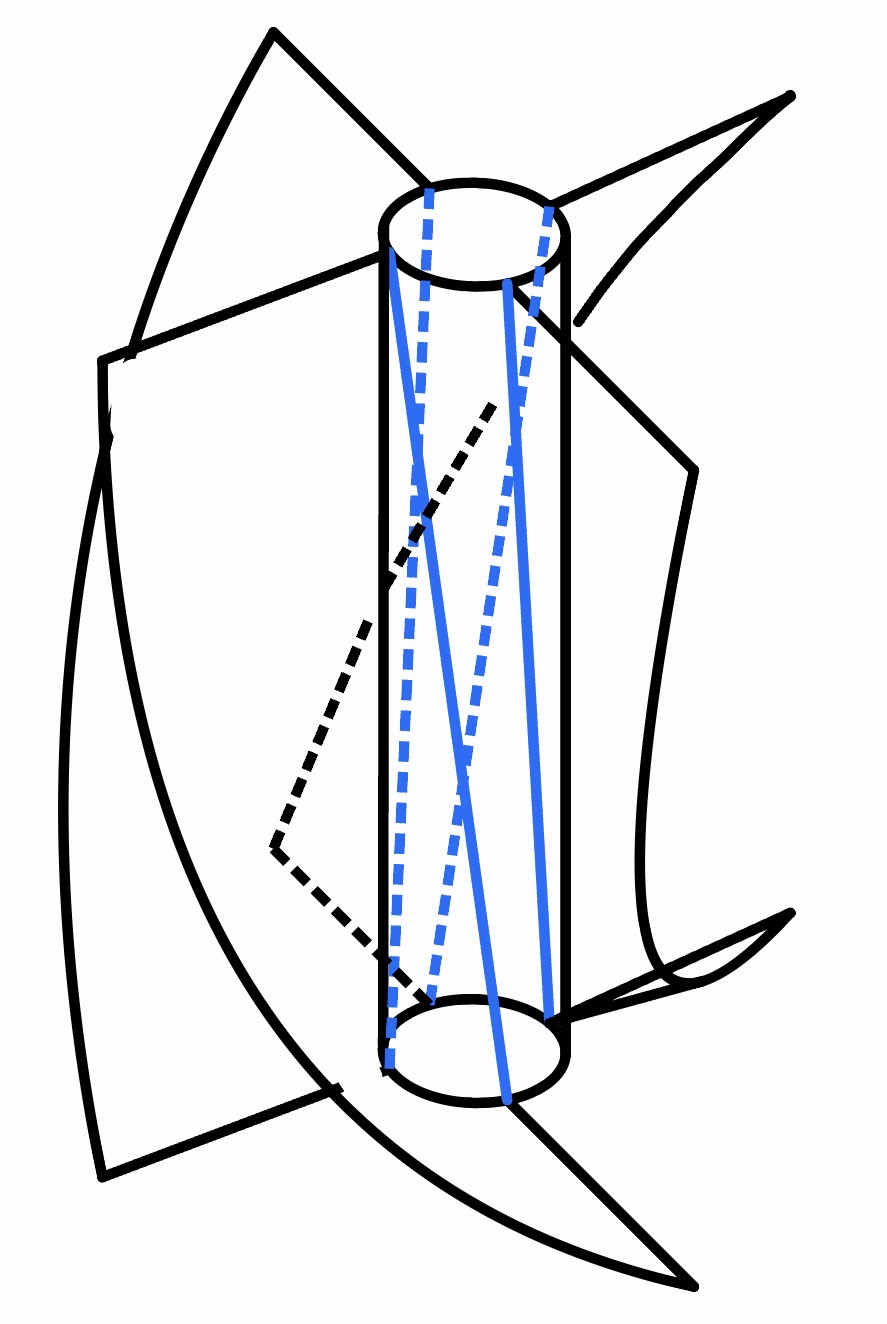}

    \caption{Left: Guts decomposition of a closed complement into the interstitial bundle $\mathcal{I}$ and the compact gut $G$. 
    Right: A local picture of the degeneracy curve in blue near a 4-pronged singular orbit}
    \label{fig: guts decomp and degeneracy}
\end{figure}

Using the topology of $M$, we can prove a stronger result about guts of the closed complement of the stable laminations $\Lambda^s$. Recall that the closed complement of an essential lamination can be uniquely decomposed up to isotopy into a union of the interstitial bundle $\mathcal{I}$ and the gut $G$. The interstitial bundle is $\pi_1$-injective and is a maximal union of I-bundles over possibly noncompact surfaces. The gut $G$ is compact and meets the interstitial bundle along properly embedded essential annuli called the interstitial annuli. See Figure \ref{fig: guts decomp and degeneracy} for an illustration.

\begin{customthm}{B}\label{thmB}
    Let $M$ be a closed, atoroidal 3-manifold. There is a triangulation $\tau$ so that there are only finitely many triangulated solid tori in $\tau$ that can arise as gut regions of stable laminations of pseudo-Anosov flows on $M$.
\end{customthm}

One immediate corollary is that we have reduced the finiteness conjecture for pseudo-Anosov flows on closed manifolds to studying Anosov flows on manifolds with boundary.

\begin{corollary}
    To prove the finiteness conjecture for pseudo-Anosov flows on closed atoroidal 3-manifolds, it suffices to show the following: For any atoroidal 3-manifold $M$ whose boundary is a collection of tori $\cup_i T_i^2$ with chosen degeneracy slopes, there are only finitely many Anosov flows on $M$ that match the chosen slopes on the boundaries.
\end{corollary}

\subsection*{Idea of the proof}
The proof of the main theorem relies on a theorem of Gabai \cite{kneser}, which states that for an atoroidal 3-manifold $M$, there is a triangulation $\tau$ so that every essential lamination and taut foliation can be isotoped to be normal to $\tau$. In other words, this fixed triangulation $\tau$ guides isotopies to a standard position. 

In the proof, we first produce the stable lamination $\Lambda^s$ of the flow by Denjoy splitting the singular foliation $\mathcal{F}^s$. Gabai's triangulation gives a combinatorial decomposition of the complementary regions of $\Lambda^s$. On the other hand, the guts decomposition gives another natural decomposition of the same complementary region. 

The main idea in this paper is to align the two decompositions via isotopy, so that the gut region of the stable lamination $\Lambda^s$ consists of finitely many combinatorial polyhedra given by the triangulation. We argue that there are only finitely many ways for leaves of the lamination to cut 3-simplices into polyhedra. Further, we argue that up to isotopy, there are only finitely many boundary interstitial annuli to these gut regions. Since degeneracy curves of pseudo-Anosov flows are isotopic to core curves of the interstitial annuli, this gives finiteness of degeneracy curves. Degeneracy curves are cables of singular orbits, so we also conclude finiteness of singular orbits. 

The proof of the finiteness of guts uses classical techniques in normal surface theory following work of Floyd--Oertel \cite{FloydOertel}, especially branched surface equations. Morally, we show that if the same interstitial annuli bound arbitrarily large solid torus guts, there is a limit object in which the Kneser branched surface carrying the stable lamination also carries a torus. This contradicts a theorem of Gabai \cite{kneser}.

\subsection*{Outline of the paper} The paper is organized as follows: In Section \ref{section: Prelim} we briefly survey pseudo-Anosov flows and normal surfaces. In Section \ref{Section: Proof}, we give a proof of Theorem \ref{main}. We give a proof of Theorem \ref{thmB} using the topology of $M$ in Section \ref{Section: Guts}.

\subsection*{Acknowledgments}
I would like to thank my advisor Dave Gabai for many insightful conversations as well as his continual support and mentorship. I also thank Ian Agol, Saul Schleimer, Thomas Barthelm\'e, Chi Cheuk Tsang and Jonathan Zung for their interest in this work and helpful conversations, and I thank Dan Margalit, Chi Cheuk Tsang, and Sam Taylor for comments on an earlier draft.

This research was partially supported by the National Science Foundation  Grant No. DMS-2304841 in Fall 2025 and Grant No. DMS-1928930, while the author was in residence at the Simons Laufer Mathematical Sciences Institute in Berkeley, California, during the semester of Spring 2026. The author learned of this problem from the Log Cabin Workshop in 2025 and would like to thank the organizers. The author is also supported by a Simons Dissertation Fellowship in Mathematics SFI-MPS-SDF-00015436.

\section{Preliminaries}
\label{section: Prelim}
\subsection{Pseudo-Anosov Flows}
\begin{figure}
    \centering
    \includegraphics[scale=0.09]{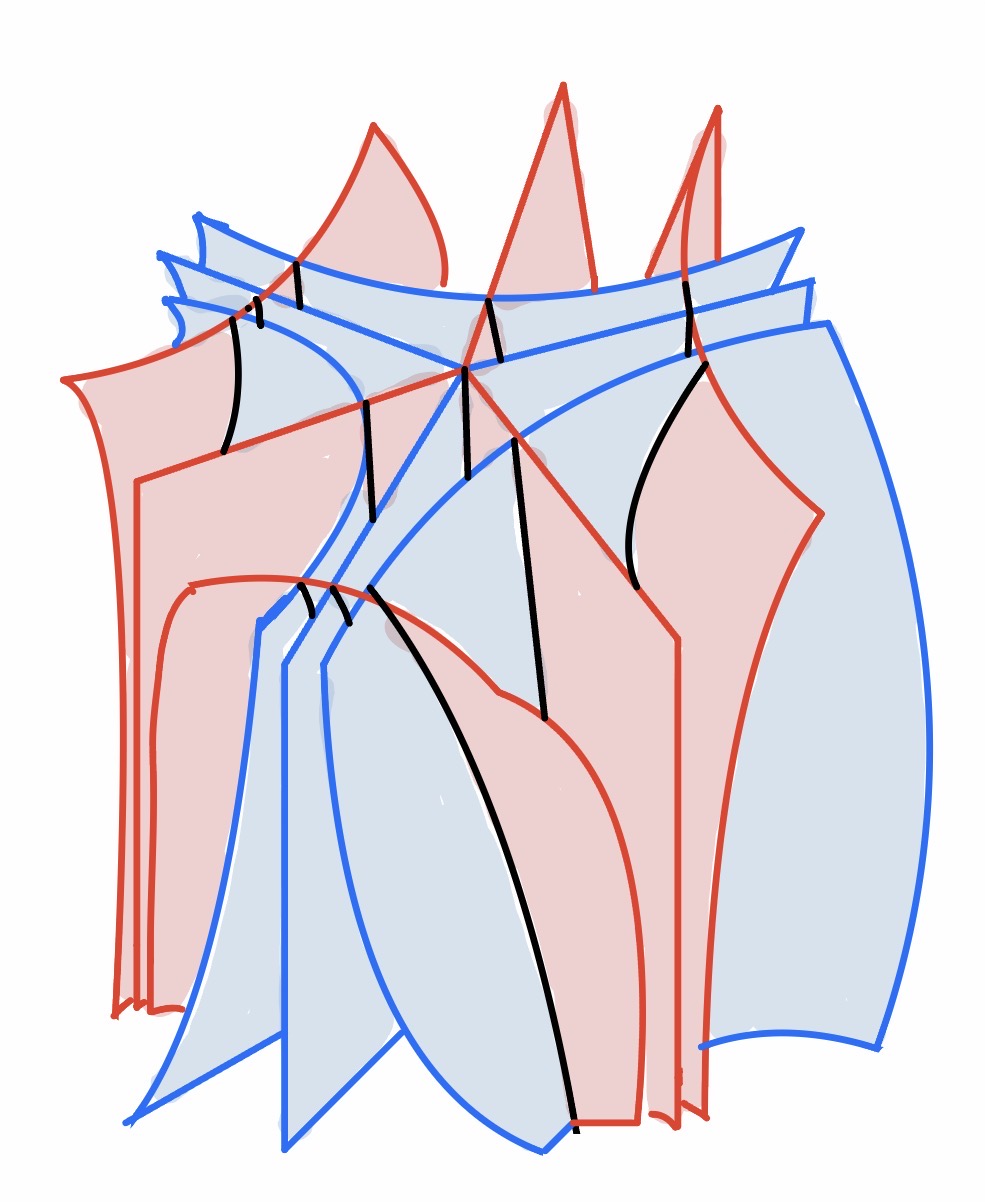}
    \includegraphics[width=0.5\linewidth]{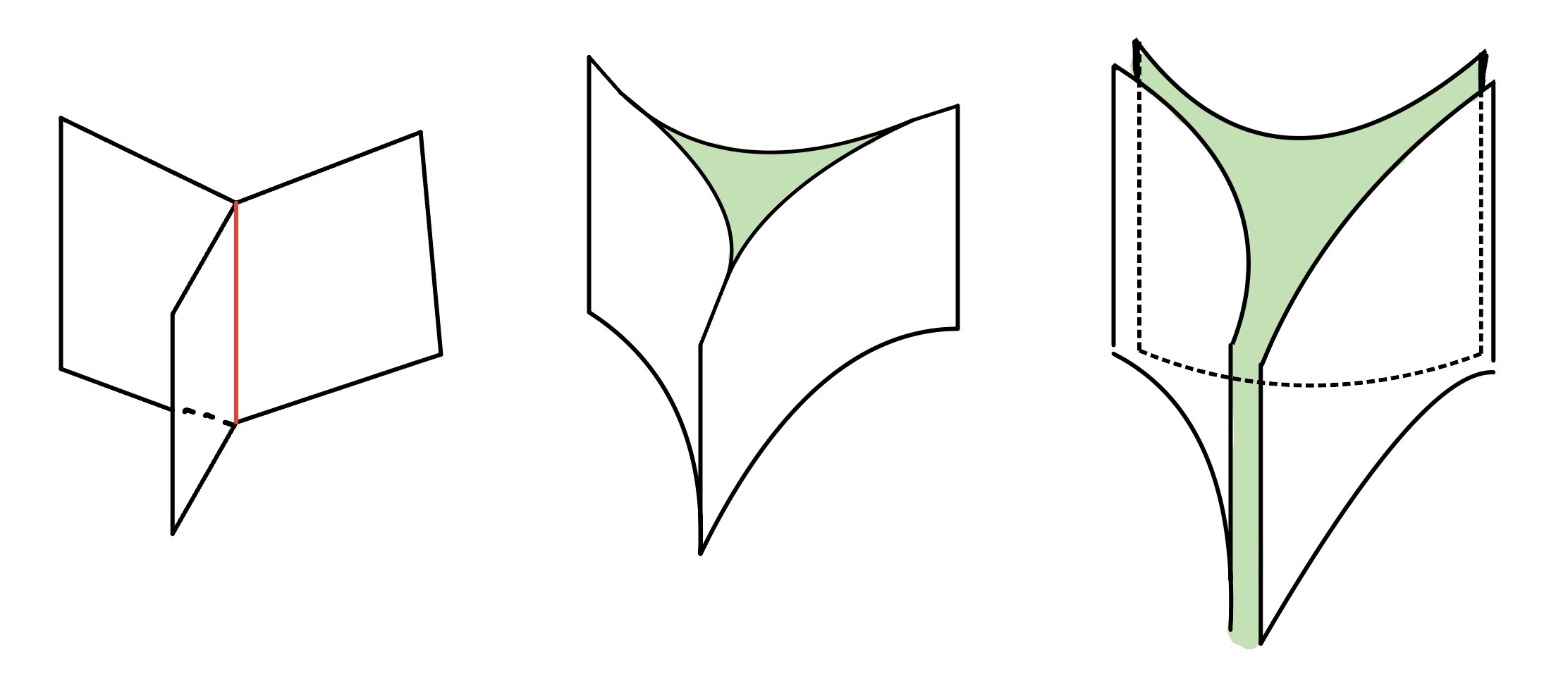}

    \caption{Left: A pseudo-Anosov flow near a singular orbit with 3-pronged singularity. Right: Denjoy splitting near a singular orbit.}
    \label{fig:prongs denjoy}
\end{figure}

We give a definition of a topological pseudo-Anosov flow following \cite{BMbook}.
\begin{definition}[Pseudo Anosov Flows]
    A \textit{pseudo-Anosov flow} is a flow $\varphi:M\times \R\rightarrow M$ satisfying the following properties:
    \begin{enumerate}
        \item There are two transverse 2D singular foliations $\mathcal{F}^s, \mathcal{F}^u$ whose leaves are saturated by orbits of $\varphi$ and intersect along orbits. They are called the \textit{stable and unstable foliations}, respectively.
        \item Singularities of $\mathcal{F}^s$ and $\mathcal{F}^u$ coincide and lie along a finite collection of periodic orbits $\alpha_i$. We call these orbits \textit{singular orbits}. A neighborhood of each singular orbit is homeomorphic to a mapping torus of a pseudo-Anosov homeomorphism of a surface near a fixed $p$-pronged singularity. See Figure \ref{fig:prongs denjoy} for an illustration of $\mathcal{F}^s$ and $\mathcal{F}^u$ near a singular orbit.
        \item In a leaf of the stable foliation $\mathcal{F}^s$, any two flow lines converge in forward time and diverge in backward time.
        \item In a leaf of the unstable foliation $\mathcal{F}^u$, any two flow lines diverge in forward time and converge in backward time.
    \end{enumerate}
\end{definition}

We consider flows up to orbit equivalence. Two flows $\varphi_1, \varphi_2$ are \textit{orbit equivalent} if there is a homeomorphism $f:M\rightarrow M$ that maps orbits of $\varphi_1$ to orbits to $\varphi_2$.

To specify the dynamics near a singular orbit, we consider also \textit{degeneracy curves}, which are the isotopy classes of multi-curves determined by the boundary of the half-leaves of $\mathcal{F}^s$ in a neighborhood of the singular orbit, as shown in Figure \ref{fig: guts decomp and degeneracy}. These multi-curves are cables of the singular orbits, and they record the rotational monodromy of the $p$-pronged singularities.

For our purposes, it is useful to consider the \textit{Denjoy splitting}: Near each singular orbit, we can blow up the stable foliation and split open the singular half-leaves to obtain a lamination $\Lambda^s$. We can similarly do this for the unstable foliation to obtain $\Lambda^u$. These laminations have no compact leaves due to the expanding dynamics of the flow, and their complementary regions are ideal polygon bundles over $S^1$. We offer a local picture in Figure \ref{fig:prongs denjoy}. For more details, see Calegari's book \cite[\S 6.6]{foliationbook}.

\subsection{Normal Surfaces} 
\begin{figure}
    \centering
    \includegraphics[scale = 0.3]{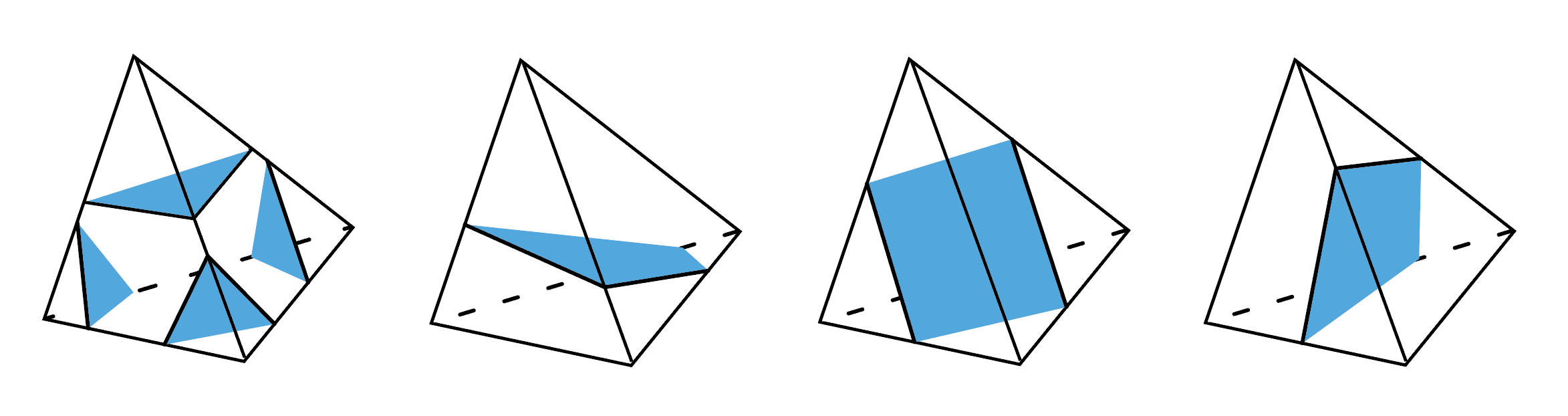}

    \caption{All 7 types of normal disks in a 3-simplex}
    \label{fig:normalSurf}
\end{figure}

The other main tool in this paper is the theory of normal surfaces. 
\begin{definition}[Normal surface]
Let $\tau$ be a triangulation of a close manifold $M$. An embedded surface $S\subset M$ is \textit{normal} if the following holds:
\begin{enumerate}
    \item It does not pass through any vertices of $\tau$, 
    \item It is transverse to the 1-complex $\tau^1$, and 
    \item For every 3-simplex $\sigma$, the intersection $S\cap\sigma$ is homeomorphic to a disjoint union of finitely many \textit{normal disks}. Normal disks arise from the intersection between affine 3-simplex and affine planes, as illustrated in Figure \ref{fig:normalSurf}. 
\end{enumerate}
Analogously, a lamination $\Lambda$ is \textit{normal} to $\tau$ if conditions (1), (2) hold, and we allow normal disks in condition (3) to intersect each 1-complex at a Cantor set instead of requiring finitely many normal disks in each 3-simplex.
\end{definition}

There are 7 types of normal disks within each 3-simplex, and we observe that no two different types of quadrilateral normal disks can coexist in the same 3-simplex. 

At its core, normality is a kind of combinatorial standard position. The Haken lemma states that every incompressible surface in an irreducible 3-manifold is isotopic to a normal surface \cite{hakenlemma}. Brittenham shows in \cite{brittenham} that if a triangulated manifold $M$ has an essential lamination, then it also has one that's normal to the triangulation. Gabai generalizes Brittenham's result in the Kneser normal form theorem, which allows us to normalize foliations and laminations to a triangulation that depends only on the ambient manifold $M$.
 We record a version of this result below.
 
\begin{theorem}[Gabai \cite{kneser}, Corollary 6.21]\label{knf}
    Let $M$ be a close orientable atoroidal manifold. There exists a triangulation $\tau$ on $M$ such that any taut foliation or nowhere dense essential lamination can be isotoped to be normal to $\tau$.
\end{theorem}

\section{Finiteness of Singular Orbits}
\label{Section: Proof}
We give a proof of Theorem A in this section. Fix an atoroidal 3-manifold $M$.

For any candidate pseudo-Anosov flow $\varphi$, we obtain its stable lamination $\Lambda^s$ by Denjoy splitting the singular stable foliation $\mathcal{F}^s$. The result $\Lambda^s$ is a nowhere dense essential lamination on $M$. To consider the stable laminations of different pseudo-Anosov flows on the same 3-manifold $M$, we invoke Gabai's Kneser normal form theorem, which produces a triangulation of $M$. 

By Gabai's result Theorem \ref{knf}, we can normalize $\Lambda^s$ with respect to a triangulation $\tau$ that depends solely on $M$ and does not depend on the choice of a flow $\varphi$. Fix this triangulation $\tau$ for the rest of this section.

Fix a pseudo-Anosov flow $\varphi$, and consider the complement of the lamination $C= M^3-\Lambda^2$. Topologically, each connected component of $C$ is an ideal polygon bundle over $\mathbb{S}^1$ by construction, created by the splitting. The boundaries of each such component of $C$ are annular leaves of $\Lambda^s$ coming from splitting singular half-leaves. 

Normalize $\Lambda^s$ to the triangulation $\tau$. The complementary region $C$ is partitioned into polyhedral components by the triangulation. Within each 3-simplex $\sigma$ of the triangulation $\tau$, there are finitely many combinatorial types of polyhedra bounded by normal disks and the boundaries $\partial\sigma$.

\begin{definition}
    Let $\sigma$ be a 3-simplex in $\tau$ and $C$ the closed complement of $\Lambda^s$. Each connected component of $\sigma\cap C$ is called a \textit{block}. The triangulation $\tau$ partitions $C$ into polyhedral blocks.
\end{definition}

\begin{figure}
    \centering
    \includegraphics[scale = 0.15]{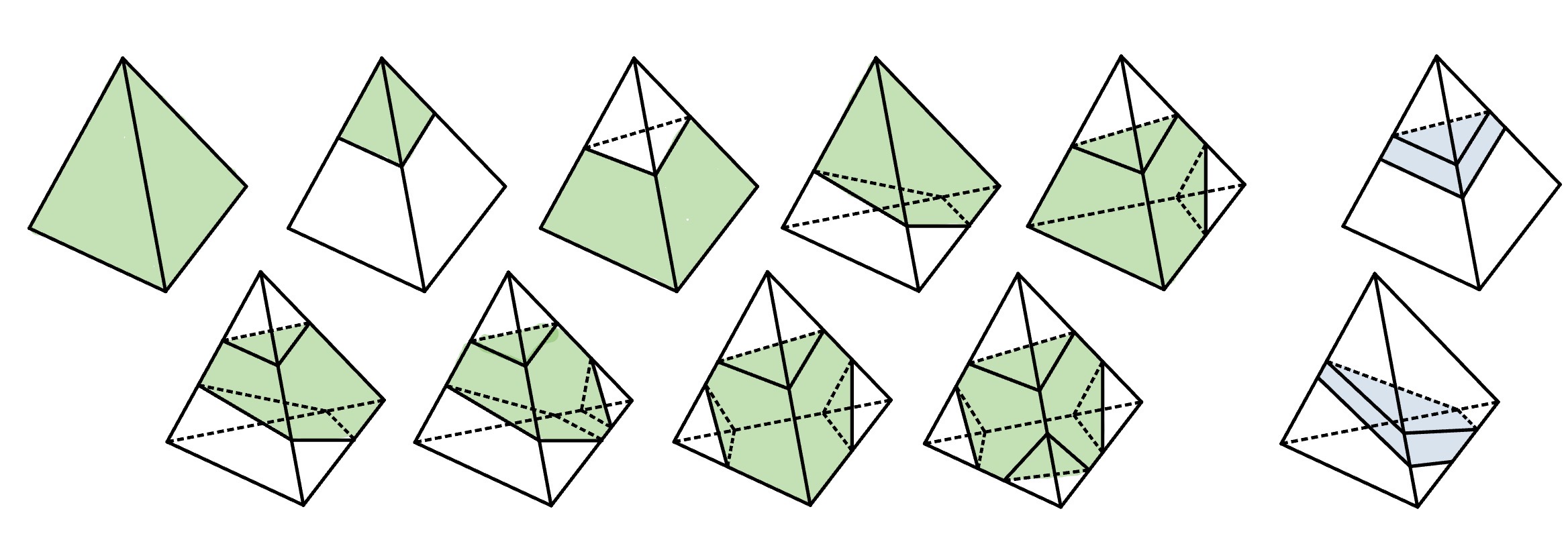}
    \caption{All possible connected polyhedral blocks  $C\cap \sigma$ in the decomposition of the normalized complementary region by the triangulation. On the left are the special blocks, and on the right are the I-bundle blocks.}
    \label{fig:blocks}
\end{figure}

Up to symmetries of 3-simplices and isotopies of normal surfaces, the diagrams in Figure \ref{fig:blocks} enumerate all possible blocks contributing to $C$. In particular, there are two types of components: I-bundles and non-I-bundle block, which we also call \textit{special blocks}. The distinguishing feature of special blocks is that for any 3-simplex $\sigma$, the complement $C\cap \sigma$ can only contain finitely many special blocks. This constitutes our first hint of finiteness.

A complementary region to the stable lamination $\Lambda^s$ decomposes into a compact gut region $G$ and some interstitial regions, which are I-bundles over non-compact surfaces  $\Sigma_j\times I$. The gut region meets the interstitial I-bundle at a collection of interstitial annuli or mobius bands. We write that a component $C_i$ of the complement of $\Lambda^s$ as $C_i=G_i\cup_{A_{ij}}\Sigma_j\times I$. This guts decomposition is unique up to isotopy, and a choice of guts decomposition amounts to choosing the interstitial annuli. 

\begin{definition}
    Let $\sigma^2$ be a 2-complex in $\tau$, and $A_{ij}$ an interstitial annulus. Each connected component of $A_{ij}\cap\sigma^2$ is called a \textit{patch}. The triangulation partitions $A_{ij}$ into quadrilateral patches.
\end{definition}

\begin{lemma}\label{chooseannuli}
There exists a choice of the interstitial annuli $A_{ij}$ such that the follow conditions are satisfied: 
\begin{itemize}
    \item $A_{ij}$ lie on the 2-skeleton of the triangulation $\tau$, 
    \item The gut $G$ inherits a partition into finitely many blocks in $C$,
    \item Every patch $P_{ijk}$ of $\sigma^2\cap A_{ij}$ borders only I-bundle blocks in the gut for each 2-simplex $\sigma^2$.
\end{itemize}
\end{lemma}
\begin{proof}
    First make an arbitrary choice for the interstitial annulus $A_{ij}$ and its corresponding gut region $G_0$. Expand $G_0$ into the interstitial I-bundle until $A_{ij}$ only intersect I-bundle blocks in the triangulation, and call the extended gut $G_1$. This is possible because there are only finitely many special blocks, so expanding the gut into the noncompact ends of the complementary region will eventually exhaust all special blocks. Define $G'$ to be the subset of $C$ defined by $G_1$ along with any I-bundle blocks whose interior intersects interstitial annuli. Each time we add in a new I-bundle block constitutes an isotopy of the interstitial annuli through that block, so the result is a valid choice of gut region and interstitial annuli. 

    Next, we argue that the gut $G_1$ consists of finitely many polyhedral components. Suppose not, then by the pigeonhole principle, there exists a 3-simplex $\sigma$ such that $\sigma\cap G_1$ has infinitely many connected components $B_k$. Choose a point $x_k$ in each of these components. Since $G_1$ is a compact subset of the complementary region $C$, there is a subsequence of $\{x_k\}$ that converges to a limit $x\in \sigma\cap G_1$. A small neighborhood of $x$ in $\sigma$ should be contained in $G_1$, but it will also contain infinitely many points of the convergent subsequence with regions not in $G_1$ separating them. This is a contradiction, and we conclude that the gut $G_1$ consists of finitely many polyhedral components from the triangulation $\tau$.

    Finally, we remark that even though guts consists of finitely many polyhedral blocks, the step to expand the gut region $G_0$ until all special blocks in the complementary region are inside the new gut $G_1$ might seem to result in arbitrarily large gut regions, measured by combinatorial complexity of how many polyhedral pieces it consists of. In other words, there isn't an obvious uniform bound on how many polyhedral blocks make up guts for stable laminations of different pseudo-Anosov flows on $M$. We show in the next section that this complexity is uniformly bounded for all pseudo-Anosov flows on $M$ if $M$ is atoroidal.
\end{proof}

For a stable lamination $\Lambda^s$, choose interstitial annuli $A_{ij}$ according to Lemma \ref{chooseannuli} for the guts decomposition of every connected component $C_i$ of the complementary region $C$. 

Recall that every interstitial annulus $A_{ij}$ is partitioned into polygonal patches $P_{ijk}$ by the triangulation, and each patch borders one polygonal block of the gut region. In the next step, we describe an inductive algorithm to simplify each gut region to a minimal position, which we define. 

\begin{definition}
    A gut region $G$ is in \textit{minimal position} if every patch of its boundary interstitial annuli borders a \textit{special block}. 
\end{definition}

\begin{algorithm}\label{alg}
    Given a gut $G$ and its interstitial annuli $A_{ij}$, run the following procedure inductively while $G$ is not yet in minimal position, one step of which is illustrated in Figure \ref{fig: one_side}.
    \begin{enumerate}
        \item Choose a patch $P$ of $A_{ij}$ that borders an I-bundle block $B$.
        \item Delete block $B$ from $G$. 
        \item Update $A_{ij}$ by an isotopy through the deleted I-bundle block.
        \end{enumerate}
\end{algorithm}
\begin{lemma}
        Given as input guts and interstitial annuli satisfying conditions in Lemma \ref{chooseannuli}, Algorithm \ref{alg} terminates in finite steps. Moreoever, at every step, we obtain valid gut regions and interstitial annuli.
\end{lemma}
\begin{proof}
    We first show that at each step, the result is still a valid gut region. In particular, we show that the algorithm only isotopes the solid torus gut and the interstitial annuli. 

    Let $P$ be the patch of the interstitial annulus chosen in step (1). It is one face of an I-bundle block $B$. By construction, $P$ is a quadrilateral with one pair of opposite sides on annular leaves of $\Lambda^s$ coming from the Denjoy splitting, and the other pair of opposite sides coincide with the 1-skeleton $\sigma^1$ of the fixed triangulation $\tau$. The interstitial annuli $A_{ij}$ is a union of these patches glued along their 1-skeleton sides.
    
    We proceed in cases. If $P$ is the only face of $B$ on the interstitial annulus, deleting $B$ results in isotoping $A_{ij}$ through $B$, a topological ball. This is illustrated in Figure \ref{fig: one_side}.
\begin{figure}

    \begin{center}
        \includegraphics[scale = 0.1]{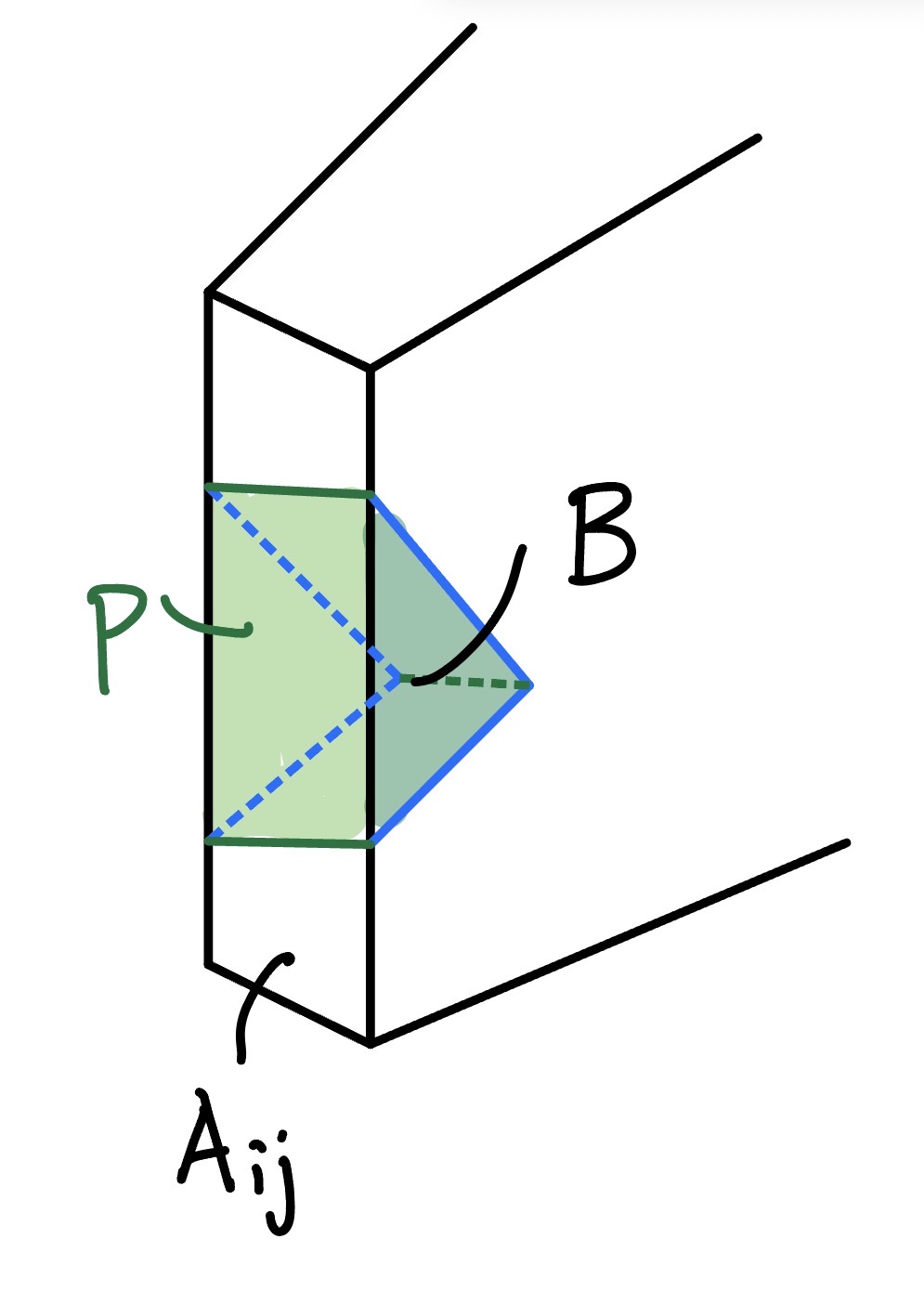}
        \includegraphics[scale = 0.1]{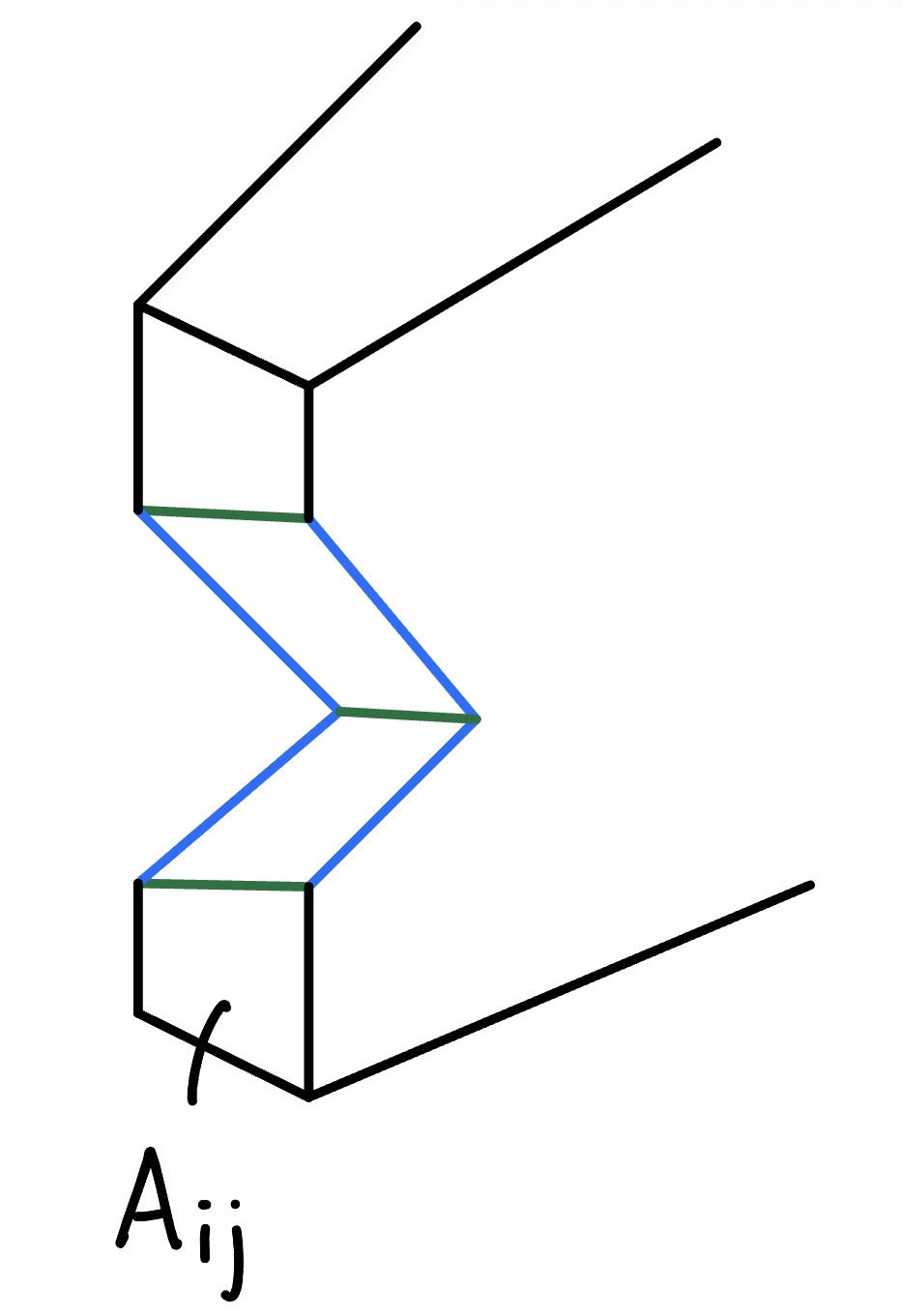}
        \caption{The result of removing block B from gut region. The interstitial annulus $A_{ij}$ is isotoped to now include the other faces of $B$ on the 2-skeleton of the triangulation.}
        \label{fig: one_side}

    \end{center}
\end{figure}
    
    If $P$ is one of two faces of $B$ on the interstitial annuli bounding the gut $G$, these two faces, call them $P, P'$, either share a side or don't. If they share a side, it reduces to the previous case. Otherwise, the two faces may be on the same annulus or different annuli. 
    

    If $P, P'$ are on the same annulus, consider the intersection $B\cap\partial G$, which is an annulus $A$ consisting of 2 patches $P, P'$ and 2 local leaves $L_1, L_2$ of the lamination. The core curve of this annulus $\gamma\subset\partial G$ bounds a disk $D$ inside the block $B\subset G$. Since $G$ is a solid torus, $\gamma$ is either nullhomotopic or the meridian on $\partial G$. If $\gamma$ is the meridian, then the disk $D$ it bounds is isotopic to a fiber of the polygon bundle over $S^1$. In particular, it informs us that the corresponding singular orbit has only 2 prongs, which is a contradiction. Therefore, $\gamma$ must not be the meridian of $\partial G$. If $\gamma$ is nullhomotopic, then removing $B$ from $G$ will separate $G$ into two connected components, one of which is a topological ball, and the other a solid torus $T$. We then isotope the interstitial annuli $A_{ij}$ through the ball component, replacing the patches of $A_{ij}$ on the ball with the patch $B\cap T$. Update $G$ to be the solid torus $T$.  

\begin{figure}

    \begin{center}
        \includegraphics[scale = 0.1]{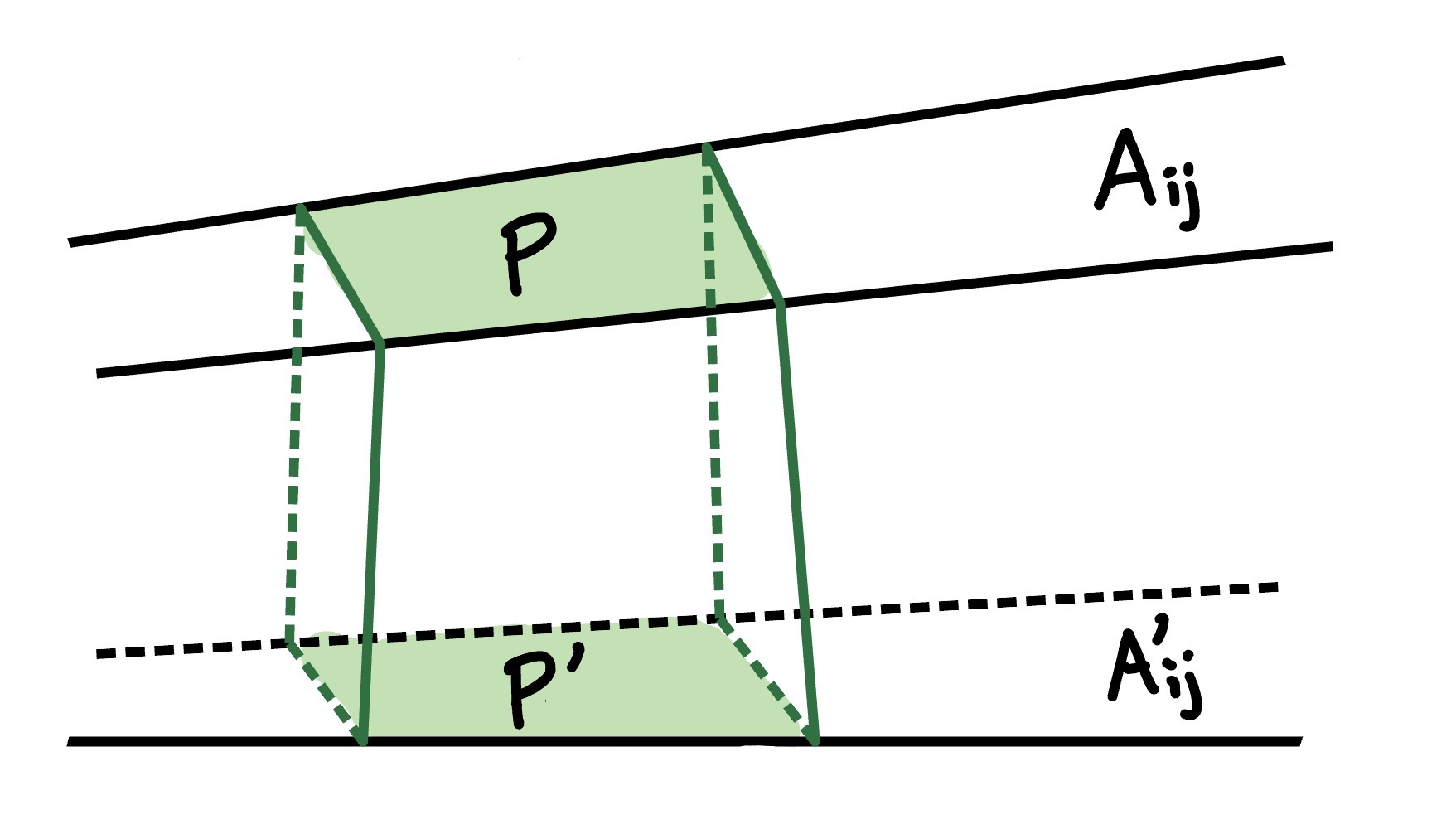}
        \includegraphics[scale = 0.1]{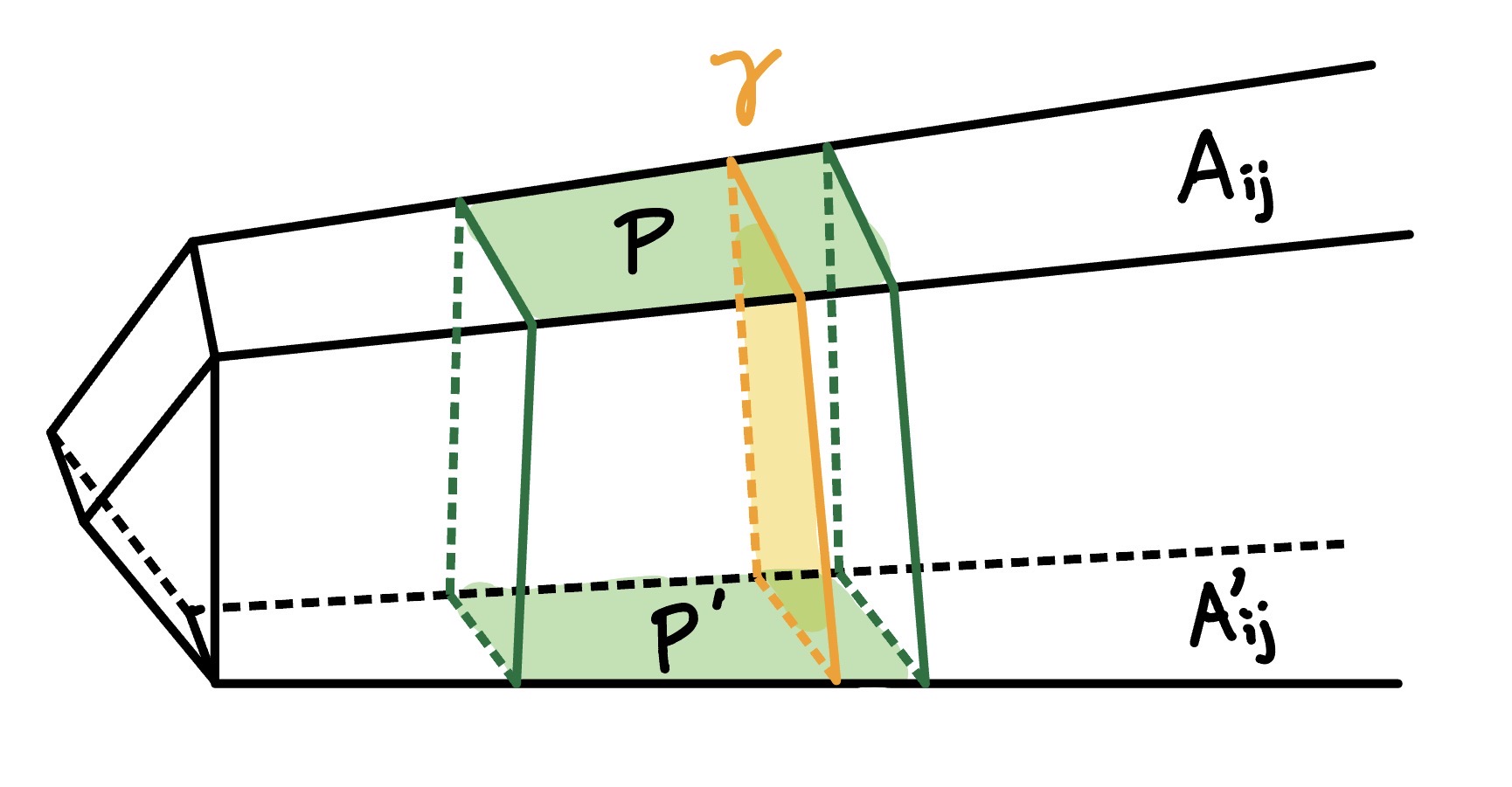}
        \caption{When two non-adjacent patches $P, P'$ belong to the same annulus, we find a topological ball to one side of $B$. Then, delete $B$ along with this ball from the guts.}
        \label{fig: two_side}

    \end{center}
\end{figure}

    Now suppose $P, P'$ lie on different annuli. Observe that $P, P'$ each has a pair of opposite sides that lie on leaves of the lamination $\Lambda^s$. Suppose $P_1$ has sides on leaves $L_1, L_2$ of $\Lambda^s$, then because $P'$ is a face of the same block $B$, $P'$ also has sides on the same leaves $L_1, L_2$ of the lamination.  But it is impossible to have two different interstitial annuli of the same gut region that meet the same leaves if the singular orbits have at least 3 prongs. Thus, $P, P'$ cannot lie on different annuli.

    If $P$ is one of 3 faces of $B$ on an interstitial annulus, all three sides must be adjacent, and we delete that block and isotope guts as before. This concludes our casework.

    By Lemma \ref{chooseannuli}, the gut $G$ consists of finitely many polyhedral blocks, and each iteration of the algorithm strictly decreases this count. Therefore, Algorithm \ref{alg} will terminate.
\end{proof}

When the algorithm terminates, every patch of the interstitial annuli $A_{ij}$ borders only a special block. We next show that there are only finitely many interstitial annuli in minimal position.

\begin{lemma}\label{finiteintersitial}
    Fixing $M$, there are finitely many isotopy classes of curves that arise as the core curves of interstitial annuli in minimal position. 
\end{lemma}

\begin{proof}
    Consider a 3-simplex $\sigma$ in the triangulation $\tau$. On each face of $\sigma$, there are at most two quadrilateral patches of any minimal interstitial annuli. This is because each patch must bound a special block, and checking against the table of special blocks, one side of a quadrilateral patch must be incident to a quadrilateral type normal disk, and the other side has only two choices for which triangular normal disk to be incident to. Therefore, there are only finitely many patches contributing to interstitial annuli, so there are finitely many combinatorial possibilities of interstitial annuli on the 2-skeleton $\tau^2$. The finiteness of core curves follow immediately.
\end{proof}

We can now prove Theorem \ref{main}.
\begin{proof}[Proof of Theorem A] Recall that we have fixed an atoroidal 3-manifold $M$ and a triangulation $\tau$ following Gabai's Kneser normal form theorem \cite{kneser}.

    Let $\varphi$ be a pseudo-Anosov flow on $M$, and let $\Lambda^s$ be its stable lamination. For every connected component $C_i$ of the complement of $\Lambda^s$, choose interstitial annuli $A_{ij}$ on the 2-skeleton of $\tau$ so that every patch borders only I-bundle blocks of the gut $G_i$. This is possible by Lemma \ref{chooseannuli}. Then use Algorithm \ref{alg} until termination. This amounts to removing  I-bundle blocks from the gut $G_i$ and isotoping the interstitial annuli $A_{ij}$ to minimal position. We remark that there may be many connected components $C_i$ of the complement of $\Lambda^s$, but after putting all interstitial annuli of every complementary region into minimal position, the interstitial annuli remain disjoint because we have only performed isotopy supported on a small neighborhood of the gut regions. By Lemma \ref{finiteintersitial}, there are only finitely many $A_{ij}$ in minimal position and finitely many core curves.
    
    Finally, observe that the core curves of the interstitial annuli $A_{ij}$ are precisely the degeneracy curves of the pseudo-Anosov flow $\varphi$, so we obtain finiteness of degeneracy curves. 
    
    Since degeneracy curves are cables of the singular orbits, we show that fixing a set of degeneracy curves $\{c_i\}$, there are only finitely many possible isotopy classes of singular orbits $\{\gamma_j\}$. Every choice of singular orbits corresponds to a collection of tori in the complement of the degeneracy curves. These tori are incompressible to the inside because of the degeneracy curves, and they are incompressible to the outside because $M$ is not a lens space. Therefore, these tori are essential tori in $M\backslash \nu(\{c_i\})$. By the JSJ decomposition, there are only finitely many such essential tori, so there are also finitely many isotopy classes of links of singular orbits.
\end{proof}

\section{Finiteness of Guts}
\label{Section: Guts}

In this section, we give a proof of Theorem B, which promotes Theorem A to a more general statement about the gut regions of the stable lamination of a pseudo-Anosov flow on a fixed 3-manifold $M$.

\begin{customthm}{B}
    Let $M$ be a closed, atoroidal 3-manifold. There are only finitely many triangulated solid tori in the triangulation $\tau$ that can arise as gut regions of stable laminations of pseudo-Anosov flows on $M$.
\end{customthm}

\begin{proof}
    Suppose not, then there are infinitely many combinatorially different gut regions $\{G_i\}$ coming from an infinite sequence of pseudo-Anosov flows $\{\varphi_i\}$ with the same interstitial annuli. In particular, the area of the torus boundary of the gut region grows unbounded. Topologically, this boundary torus $T^2$ has 1 or more parallel interstitial annuli.

     Using a result of Gabai \cite{kneser}, there is a finite collection of essential branched surfaces carrying all the weak stable laminations $\Lambda^s_i$ . By the pigeonhole principle, there is an essential branched surface $\Sigma$ that carries infinitely many stable laminations in the sequence $\{\Lambda^s_i\}$. We restrict to that subsequence.
     
     For any closed surface $S$ carried by $\Sigma$, we can write a system of linear equations $Aw=0$, where the matrix $A$ records how the sectors of the branched surface interact at the branch locus, and the vector $w$ consists of weights $w_k\in\Z$ representing the number of local sheets carried by a sector of $\Sigma$. When $S$ has boundary along the branch locus, we write $Aw=b$, where an entry $b_k$ records that along the k-th edge of the branch locus, we see $b_k$ pieces of the boundary.

    Let $\{T_i, \gamma\}$ be the sequence of boundary tori of gut regions that increase in area despite having the same interstitial annuli $\gamma$. Let $S_i = T_i\setminus\gamma$ be the part of the torus $T_i$ that coincide with the annular leaves of the stable lamination coming from splitting singular leaves. Computing the Euler characteristic, we have the following. 
    \[\chi (S_i) = \chi(T_i)-\chi(\gamma)=\chi(T_i)-0=0\]

    Consider the branched surface equations for $S_i$. Since $S_i$ is not a closed surface, we can write its branched surface equations $Aw^{(i)} = b^{(i)}$. Because the interstitial annuli are the same for all $T_i$, $b^{(i)} = b$ for some fixed $b$ for all $i$. This is because the vector $b^{(i)}$ only records boundary information for $S_i$. 

    Next, we normalize the weight vectors
    \[x^{(i)}=\frac{w^{(i)}}{\norm{w^{(i)}}},\]
    and pass to a convergent subsequence in the normalized sequence $x^{(i)}\rightarrow x$.
    Since the area grows unbounded, the norm $\norm{w^{(i)}}$ goes to infinity in the limit, so that \[Ax^{(i)}=\frac{b}{\norm{w^{(i)}}}\rightarrow 0,\]
    and the limit satisfies 
    \[Ax=0.\]

    Work of Floyd-Oertel \cite{FloydOertel} tells us that a non-negative solution $x$ to the branch equations implies that $x$ defines a measured lamination $\Lambda$ carried by $\Sigma$. Since the Euler characteristic is linear in the weights $w^{(i)}$, $\chi(\Lambda)=0$. 

    The limit $x$ lies in the solution space $\{w\geq 0: Aw=0\}$, where every solution is a positive linear combination of extremal rays $r_j$ in the solution space (similar to fundamental solutions): $x=\sum_j c_j r_j$, where $c_j\geq 0$. Since every $r_j$ is a non-negative solution to the branch equations, it defines a measured lamination $\Lambda_j$. Because Euler characteristic is additive, and $\Sigma$ carries no sphere leaves of $\Lambda_j$ that might contribute positive Euler characteristic, we have $\chi(\Lambda_j)=0$. But we know from linear algebra that the extremal rays consists of rational entries because $A$ has integral entries, so $r_j$ actually define embedded surfaces. By the Euler characteristic, $r_j$ defines an embedded torus. But this is a contradiction, because the branched surface $\Sigma$ does not carry tori.  

    Therefore, we conclude that this infinite sequence of guts can not exist in the first place. There are only finitely many combinatorial possibilities for triangulated guts of pseudo-Anosov flow on atoroidal $M$.
\end{proof}

\begin{corollary}\label{lengthbound}
    Let $M$ be a closed atoroidal 3-manifold. There is a uniform bound on the length of the singular orbits of pseudo-Anosov flows on $M$.
\end{corollary}
We offer two proofs of this corollary. The first one is an observation of Saul Schleimer \cite{Saulcomment}.

\begin{proof}[First Proof]

Because the gut regions are normalized and thereby triangulated by $\tau$, we can invoke a result of Lackenby \cite{lackenby}, which gives a uniform bound on the length of the core curve of a triangulated solid torus. This gives the desired corollary since singular orbits are the core curves of triangulated solid tori guts. 
\end{proof}
Another viewpoint of Corollary \ref{lengthbound} is through the degeneracy curves of singular orbits.
\begin{proof}[Second Proof]
 By Theorem \ref{main}, degeneracy curves of singular orbits are finite combinatorial objects on the 2-skeleton $\tau^2$ of the triangulation passing each 2-complex at most twice, so they also are uniformly bounded in length. Degeneracy curves are cables of the singular orbits, so the singular orbits have uniformly bounded length. 
\end{proof}

\section{Further Questions}

One direction of future investigation is obtaining an effective bound on the number of isotopy classes of links that arise as singular orbits of pseudo-Anosov flows for a manifold $M$. Let the set of candidates for singular orbits determined by Theorem \ref{main} be denoted $S(M)$.

More precisely, here are some questions we would need to answer to this end: 

\begin{question}\label{question 1}
    Given an atoroidal 3-manifold $M$, how do we characterize the triangulation given by Gabai's Kneser normal form theorem? How many tetrahedra is sufficient? How do we algorithmically find such a triangulation? 
\end{question}

One motivation for Question $\ref{question 1}$ is that the proof of Theorem \ref{main} is constructive given the triangulation $\tau$ from Gabai's Kneser normal form result \cite{kneser}. However, Gabai's proof of the existence of $\tau$ as in Theorem \ref{knf} uses a compactness argument to prove that every essential lamination on a closed atoroidal 3-manifold $M$ is carried by one of finitely many essential branched surfaces. It is therefore not constructive. Agol--Li \cite{Agol_2003} gives an algorithmic rendering of the finiteness of carrying essential branched surfaces using 1-efficient triangulations. However, there is not an effective bound on the combinatorial complexity of these branched surfaces.

On the other hand, Gabai shows in the same paper \cite{kneser} that given an arbitrary triangulation $\Delta$ on a closed 3-manifold $M$, there's an isotopy to normalize a lamination $\Lambda$ to $\Delta$ after performing a bounded number of evacuations, which he defines. It's plausible that we can enumerate the possible singular orbits by choosing a preferred triangulation apart from $\tau$ and studying the possible evacuations. 

\begin{question}
    Is there an efficient algorithm to generate all possible singular orbits of pseudo-Anosov flows on $M$ given an arbitrary triangulation of $M$? 
\end{question}

\begin{question}
    Given a subset $A\subset S(M)$, is there an efficient algorithm to determine whether $A$ could arise as the link of singular orbits of a pseudo-Anosov flow?  
\end{question}

Recent work of Schmalian \cite{SchmalianObstruction} gives first examples of cusped hyperbolic manifolds that are not the non-singular part of a pseudo-Anosov flows, giving an obstruction using the Heegaard-Floer homology of the double branched cover. This obstructions seems effective for small manifolds: Schmalian shows that among the first 100 manifolds in the orientable cusped manifold census, 19 have no veering triangulations, compared to 20 appearing in the veering census. There's work in progress of Hall \cite{Hall} exhibiting large classes of manifolds that do not admit veering triangulations or pseudo-Anosov flows without perfect fits.

\bibliography{references}
\bibliographystyle{alpha}

\end{document}